\documentclass{article}
\usepackage{graphicx}
\usepackage{amssymb}
\newtheorem{theo}{Theorem}[section]
\newtheorem{cor}[theo]{Corollary}
\newtheorem{lem}[theo]{Lemma}
\newtheorem{prop}[theo]{Proposition}
\newtheorem{defi}[theo]{Definition}

\newcommand{\proof}{\noindent{\bf Proof}\par\nobreak}

\renewcommand{\>}{\rangle}

\newcommand{\qed}{\ \hfill$\square$\smallskip}

\title{Subgroup theorem for valuated groups and the CSA property}
\author{Abderezak OULD HOUCINE}
\date{}
\begin{document}
\maketitle
\begin{abstract}
A   valuated group with normal forms is a group  with an
integer-valued length function satisfying some Lyndon's axioms
\cite{Lyndon-length} and an additional axiom considered by Hurley
\cite{hur2}. We prove a  subgroup theorem for valuated groups with normal forms analogous to   Grushko-Neumann's theorem. We study also the CSA property in such groups.
\end{abstract}
\maketitle
\section{Introduction}
Let $(\Lambda, \leq )$ be an (totaly) ordered abelian group, and
$G$  a group with a length function $\ell :G\rightarrow \Lambda $.
For $x,y \in G$, we  let $c(x,y)=\frac{1}{2}(\ell(x) +\ell(y)
-\ell( xy^{-1}) )$. We notice that we may have  $c(x,y) \not \in \Lambda$, but we may assume that we are working in the divisible ordered abelian closure of $\Lambda$ (see \cite{chiswell} for more details).

We say that $\ell$ is a \emph{Lyndon length function}, if it
satisfies the following axioms considered by Lyndon \cite{Lyndon-length}:\\

\noindent $\mathcal{A}_{1}$. $\ell( 1) =0$,

\noindent $\mathcal{A}_{2}$. for all $x\in G$, $\ell( x^{-1})
=\ell(x),$

\noindent $\mathcal{A}_{3}$. for all $x,y,z \in G$, $c(x,y) \geq
\mbox{min}
\{c(x,z),c(z,y)\}$,\\

\noindent and in that case $(G,\ell)$ is called a
\emph{$\Lambda$-valuated group}. If $\Lambda$ is $\mathbb{Z}$ with
the usual ordering, we call $(G,\ell)$  a \emph{valuated group}. We shall use the notation $\ell$ for  length functions unless otherwise indicated.

In \cite{chis2}, Chiswell showed that a valuated group $(G, \ell)$,  assuming that $c(x,y)$ is always an integer,  acts on a tree $T$ in such a way that  $\ell(g)$ is the tree distance between $p$ and $gp$ for some suitable vertex $p$ of $T$. Conversely, if $T$ is a tree and $p$ is vertex of $T$, and $G$ is a group acting on $T$, then $(G, \ell)$ is a valuated group, $\ell(g)$ being the tree distance between $p$ and $gp$. Hence, the subject of valuated groups fits in the theory of groups acting on trees.

We are concerned in this paper with a restricted class of valuated groups. We are interested in valuated  groups $G$,
which satisfy the following additional axiom considered by Hurley
\cite{hur2}:
$$\mathcal{A}_4.~G \hbox{ is generated by the set } \{x \in G~|~\ell(x) \leq 1\}.$$
\indent A valuated group $G$ satisfying $\mathcal{A}_4$ is
called a \emph{valuated group with normal forms}.  Free groups, free products with amalgamation and HNN-extensions are the typical examples of valuated groups with normal forms. Proposition \ref{prop-Hurley-normalform}
below, shows that every element of a valuated groups with normal forms  has normal
forms, having several properties comparable to those of free
product with amalgamations and HNN-extensions.

It is very pleasant  to work directly in the class of valuated
groups with normal forms for several raisons. For instance,
because this class contains free groups, free product with
amalgamation and HNN-extensions, and  in several times the same
proof works for all such groups as it depends generally on some
properties of normal forms. Therefore, this  provides a unifying framework in which we can study both free products with amalgamation and HNN-extensions.

Lyndon \cite{Lyndon-length} has introduced   integer-length functions on groups, satisfying  some
 axioms including $\mathcal{A}_1, \mathcal{A}_2, \mathcal{A}_3$, to axiomatize the argument of  Nielsen's proof
  of the subgroup theorem for free groups. He has proved a splitting  theorem (Theorem
\ref{theo-decomp-lyondon} below), which gives a new
 proof of the Nielsen subgroup theorem, the Kurosh subgroup theorem, and which gives more information about the
 restriction of the natural length on subgroups of free products.   Hurley \cite{hur2} has studied groups with normal forms (called NFS-groups), and has shown that such groups can be
obtained by considering Lyndon's axioms and  the additional axiom
$\mathcal{A}_4$ above.

By introducing a new axiom ($\mathcal{A}_5^*$ below), we prove a subgroup theorem for valuated groups analogous to
Grushko-Neumann's theorem (Theorem \ref{princip-theo}).  This gives at least an uniform statement and an uniform proof of
theorems about subgroups in free groups, free product with
amalgamation and HNN-entensions. We prove a splitting theorem for valauted groups satisfying a special axiom ($\mathcal{A}_0^*$ below), and we study centralizers in valuated groups with normal forms. 

We will be also  interested in this paper, with the CSA property. A
subgroup $H$ of a group $G$ is {\em conjugately separated} in $G$,
or {\em malnormal} in $G$, if $H \cap H^x =1$ for every $x
\in{G\setminus H}$. A {\em CSA-group} (``Conjugately Separated
Abelian'') is a group in which every maximal abelian subgroup is
malnormal.  It is known  that a
CSA-group with an {\em involution}, i.e. an element of order $2$,
must be abelian. Following the notation of
\cite{MyasnikovRemeslenikov96}, we denote by {\em CSA$^*$-group}
any CSA-group without involutions. The class of CSA-groups contains free groups,  and more generally torsion-free hyperbolic groups
\cite{Gromov87}, groups acting freely on $\Lambda$-trees
\cite{Bass91, chiswell}, and limit groups
\cite{GaglioneSpellman93, Remeslennikov89}.

In \cite{GildKharlamMyas95, Eric-Houcine},   conditions on free product with amalgamation and
HNN-extension to be a CSA$^{*}$-groups are given. We generalize in this paper that results to valuated groups with normal forms.

If $(G, \ell)$ is $\Lambda$-valuated group, it follows from the axioms, and we leave the proof to the reader, the following
properties (which will be used freely without explicit reference):\smallskip

\noindent $(i)$ $c(x,y)=c(y,x)$ and $\ell(x) \geq 0$, for all $x,y
\in G$.

\noindent $(ii)$ $| \ell(x) -\ell(y) | \leq \ell(xy) \leq
\ell(x)+\ell( y) $, for all $x,y,z \in G$, where $|.|$ is the
absolute value in $\Lambda$. In particular, if $\ell(y)=0$ then $\ell(xy)=\ell(x)$.

\smallskip
We define  $B=\{x \in G ~|~\ell(x)=0\}$. Then it follows from
$\mathcal{A}_{1}$-$\mathcal{A}_{2}$ and $(ii)$ that $B$ is a subgroup of $G$.

We consider also the following  type of length functions.
Let $G$ be a group generated by some set $S$. Then one defines a
length function $\ell_S$ over $G$, called \emph{word
length}, as follows.  For every $g \in G$, we let $\ell_S(g)$ to be the
smallest length of words $w$ over $S^{\pm 1}$ such that $g=w$.
Then $\ell_S$ takes its values in $\mathbb{N}$, and one can check
easily that $\ell_S$ satisfies $\mathcal{A}_1$, $\mathcal{A}_2$
and $ \ell_S(xy) \leq \ell_S(x)+\ell_S( y) $, for all $x,y \in G$.

The paper is organized as follows.  In the next section we gives some examples of valuated groups. Section 3 is devoted to preliminaries on valuated groups. In Section 4, we show a subgroup theorem for valuated groups with normal forms. We prove in Section 5, a splitting theorem for valuated groups satisfying axiom $\mathcal{A}_0^*$; which is a weak verion of axiom $\mathcal{A}_0$. Section 6 concerns conjugacy in valuated groups with normal forms, and we study in Section 7 centralizers in such groups. We prove in Section 8, a theorem which gives sufficients conditions for a valuated group to be a CSA$^*$-group.

We denote by $X^\#$ the set of notrivial elements of $X$, and $g^h$ denotes $h^{-1}gh$. 
\section{Examples  \& other axioms}

We give now some examples of $\Lambda$-valuated groups.

\smallskip \noindent \textbf{(1)}  Typical examples.

\noindent $\bullet$ \emph{Free groups}. The group $(F,\ell_X)$,
where $F$ is a free group with basis $X$ and $\ell_X$ is   the word
length function, is  a valuated group satisfying the following axiom
(introduced by Lyndon):
$$\mathcal{A}_0. ~~\ell(x^2) > \ell(x),  \hbox{ for every } x \not =1.$$
Conversely, Lyndon \cite{Lyndon-length} has shown that a valuated group
satisfying $\mathcal{A}_{0}$ is also a free group. Remark that in
this case we have $B=\{1\}$.

\smallskip
\noindent $\bullet$\emph{ Free products}. Let $G = G_1*G_2$. Put
$S= (G_1 \cup G_2)\setminus \{1\}$. Then $(G, \ell_S)$ is a
valuated group which satisfies the following axiom (introduced by
Lyndon):

$$\mathcal{A}_5. ~~c(x,y)+c(x^{-1},y^{-1}) > \ell(x)=\ell(y) \Rightarrow x=y$$

Of course  $\ell_S(g)$ is the length of normal forms of $g$. Remark also that in this case we have $B=\{1\}$.

\smallskip
\noindent $\bullet$ \emph{Amalgamated free product}. Let $G =
G_1*_AG_2$. Let $\ell(g)$ be the length of the  normal form of $g$
if  $g \not \in A$ and  $\ell(g)=0$ if $g \in A$. Then $(G, \ell)$
is a valuated group satisfying the Chiswell's axioms
\cite{chiswell3}:

$$C1'. \hbox{ If } \ell(x) \hbox{ is even and } \ell(x) \neq 0, \hbox{ then } \ell(x^2)> \ell(x),$$
$$C2. \hbox{ for no }x \hbox{ is } \ell(x^2)=1+\ell(x). $$

Conversely, Chiswell \cite{chiswell3} has shown that a valuated
group with normal forms $(G, \ell)$ satisfying the above axioms is
a free product of a family $\{G_i | i \in I\}$ with a subgroup $A$
amalgamated, such that $\ell$ is the natural length function
relative to this decomposition. Note that we have $B=A$.

\smallskip \noindent $\bullet$ \emph{HNN-extensions}. Let
$G^*=\<G,t | A^t=B\>$ be an HNN-extension. Let $\ell(x)$ be the
number of occurrences of  $t^{\pm}$ in the normal form of $x$.
Then $(G, \ell)$ is a valuated group with normal forms, with in this case $B=G$.
Remark that we have  $\ell(t)=1$ and $\ell(t^2)=2$, unlike the
the case of amalgamated free products where we have $\ell(x^2)\leq
1$ for every $x$ such that $\ell(x) \leq 1$.

\smallskip
\noindent\textbf{(2)} Model of the universal theory of
non-abelian free groups. 

Let  $\cal M$ be a model of the universal
theory of nonabelian free groups. Then $\cal M$ embeds in an
ultrapower $^*F$ of  $F_2$ which   can be equipped
 with a Lyndon function $^*\ell$   taking its values in
an ultrapower $^*\mathbb{Z}$ of $\mathbb Z$ and satisfying $\mathcal{A}_0$. We get that
$(\mathcal{M}, {^*\ell |\mathcal{M}})$ is a
$^*\mathbb{Z}$-valuated group satisfying $\mathcal{A}_0$.

This viewpoint  was used by  Chiswell and Remeslennikov
\cite{chis-rem},  to give a new proof of a theorem of Appel and
Lorents,  about the solutions  of equations with one variable in free
groups.

\smallskip
\noindent\textbf{(3)} Groupes acting on
$\Lambda$-tree.

In \cite{chis2, chiswell}  Chiswell
considered valuated groups and has shown  that such groups can be
obtained from their action on a suitable tree. More generaly he
has shown that if  $(G, \ell)$ is $\Lambda$-valuated group, then $G$ acts by isometry on $\Lambda$-tree and if  $(X, d)$ is a $\Lambda$-tree and if $G$ is a
group acting by isometry over $X$, by defining for  $x \in X$
$\ell(g)=d(x,gx)$, then $(G,\ell)$ is a $\Lambda$-valuated group.

\smallskip
\noindent\textbf{(4)} Some linear groups.

Let $(K, v)$ be a fields with a discrete
valuation $v$. Then, by a result of Serre \cite{serre-tree},   GL$_2(K)$ acts on tree. Therefore,  it can be
equipped with a  Lyndon length function. Chiswell \cite{chis4} has given the
explicit form of that function,  and has showen that the corresponding tree is the same as that constructed by Serre. If we take   $\mathbb{Q}$ with the
$p$-adic valuation $v_p$, we get  GL$_2(\mathbb{Q})$ as a valuated
group.

\section{Preliminary}
\subsection{Valuated groups}

We define some interesting subsets and relations of a
$\Lambda$-valuated group $(G, \ell)$.  We let $$N_G=\{g \in G~ |~
\ell(g^2) \leq \ell(g)\},$$ and we denote it $N$ if there is no risque of ambiguity. We let $\equiv$ to be the following
relation, defined in $N$, by:
 $$x\equiv y~
\hbox{ if and only if } ~ \ell(x^{-1}y) \leq \ell( x) =\ell( y)
.$$ \indent This relation is  du  to Lyndon \cite{Lyndon-length}. One
can check easily that $\equiv $ is an equivalence relation on $N$.
Indeed, obviously  $\equiv$ is reflexive and symmetric. Let
$x,y,z \in N$ such that $x \equiv y$,  $y \equiv z$. Then
$\ell(x)=\ell(y)=\ell(z)$ and
$$c(x,y)=\hbox{$\frac{1}{2}$}(\ell(x) +\ell(y) -\ell( x^{-1}y) ) \geq \hbox{$\frac{1}{2}$}\ell(x),$$
$$c(y,z)=\hbox{$\frac{1}{2}$}(\ell(y) +\ell(z) -\ell( y^{-1}z) ) \geq \hbox{$\frac{1}{2}$}\ell(y),$$
and by axiom $\mathcal{A}_3$ we have $c(x,z)
\geq\frac{1}{2}\ell(y)$. Hence $\ell(x^{-1}z) \leq \ell( x) =\ell(
y)=\ell(z)$, thus $x \equiv z$.

We denote by $N^*(x)$ the equivalence class of $x$ under $\equiv$.
Then we see that $N^*(1)=B$. We let $N(x)=N^*(x) \cup \{1\}$.

Let $U$ be a subset of $G$. Let $(u_{1},\dots ,u_{n})$ be a
sequence in $ U^{\pm 1}$. We say that  $(u_{1},\dots,u_{n})$ is
\emph{pseudo-reduced} if:

\smallskip
$(i)$ $u_{i}u_{i+1}\neq 1$ , $u_{i}\neq 1,$

$(ii)$ if $u_{i},u_{i+1}\in  U^{\pm 1} \cap N$, then $
u_{i}\not \equiv u_{i+1}$.

\smallskip
We need in the sequel the following theorem of Lyndon, which can be
extracted from its results in \cite{Lyndon-length}.

\begin{theo} \emph{\cite{Lyndon-length}}\label{theo-decomp-lyondon}
Let $(G,\ell )$ be a valuated group satisfying
$$
\mathcal{A}_1^*. ~\ell(x)=0 \Rightarrow x=1,
$$
$$
\mathcal{A}_5. ~~c(x,y)+c(x^{-1},y^{-1}) > \ell(x)=\ell(y)
\Rightarrow x=y.
$$
Then

\emph{(1)} For every $x \in N,$ the set  $N(x)$ is a
subgroup of $G$.

\emph{(2)} There exists a generating set  $U$ of $G$ such that:

$(i)$ for every pseudo-reduced sequence $(u_{1},\dots,u_{n})$ of
$U^{\pm 1}$ we have: $$ \ell(u_{1}\cdots u_{n}) =\sum _{i=1}^{n}\ell
(u_{i}) -2\sum _{i=1}^{n-1}c(u_{i},u_{i+1}^{-1}),$$

$(ii)$ we have  $ G=F*_{i}G_{i}$, where $F$ is a free group having a
basis $X \subseteq U$ and  $G_{i}=N(x)$ for some $x\in
U\cap N$. \qed
\end{theo}

We end this section with  the following lemma,  needed in the sequel,
and which generalizes  the one proved by Lyndon \cite{Lyndon-length}
in case of valuated groups satisfying $\mathcal{A}_1^* ,
\mathcal{A}_5$. The proof is an exact copy of the one in
\cite{AHoucine}, and it is left to the reader.

\begin{lem}\label{lem-suite}
Let $(G, \ell )$ be a $\Lambda$-valuated group. Let
$(g_{1},\dots,g_{n})$, $ n \geq 2$, be a sequence in $G$
satisfying:$$ c(g_{i-1},g_{i}^{-1})+c(g_{i},g_{i+1}^{-1})< \ell
(g_{i}),  \qquad 1<i<n. $$

Then $\ell( g_{1}\cdots g_{n}) =\sum _{i=1}^{n}\ell( g_{i}) -2\sum
_{i=1}^{n-1}c(g_{i},g_{i+1}^{-1}). $ \qed
\end{lem}

\subsection{Valuated groups with normal forms}

Hurley \cite{hur1}, after studying some groups with
normal forms, have noticed that the  typical examples cited
earlier and the groups that he have studied satisfy the axiom
$$\mathcal{A}_4.~G \hbox{ is generated by the set } \{x \in G~|~\ell(x) \leq 1\}$$

As it was said in the introduction, we call  a valuated group
$(G,\ell)$, a \emph{valuated group with normal forms} if
$(G,\ell)$ satisfies $\mathcal{A}_4$.

We let $S=\left\{ x \in G~|~\ell(x) \leq 1\right\}$. A sequence
$(s_1,s_2, \dots,s_n)$ of $S$ is said to be \emph{$S$-reduced}
if  $s_i .s_{i+1} \not \in S$ for every $1\leq i\leq n-1$.

The following proposition shows that every element in a valuated group with normal forms has normal forms. His proof can be extracted from \cite{hur2}, but for completeness we provide a proof.

\begin{prop}\emph{\cite{hur2}}\label{prop-Hurley-normalform} Let $(G, \ell )$ be a
valuated group with normal forms. Then:

\emph{(1)} If  $(s_{1},\dots,s_{n})$,  with $n\geq 2$, is a
$S$-reduced sequence,  then  $$\ell( s_{1} \cdots s_{n}) =\ell
_{S}(s_{1} \cdots s_{n})=n.$$

\emph{(2)} For every  $g\in G \setminus B$, there exists a
$S$-reduced sequence   $(s_{1},\dots,s_{n})$  such that $g=s_{1}
\cdots s_{n}$ and $\ell( g) =\ell _{S}(g)=n$. Thus if  $(s'_{1},\dots,s'_{m})$  is a  $S$-reduced sequence  such
that  $g=s'_{1}\cdots s'_{m}$, then $m=n$.
\end{prop}

\proof 

(1) Since $(s_1, \cdots, s_n)$ is $S$-reduced, we get $c(s_i, s_{i+1}^{-1})=0$ for $1\leq i \leq n-1$. Therefore, by Lemma \ref{lem-suite}, $\ell(s_1 \cdots s_n)=\sum_{i=1}^{i=n}\ell(s_i)$.

(2) The existence of the $S$-reduced sequence  follows from the fact that $G$ is generated by $S$ and the rest follows from (1). \qed

\begin{defi}
 A \emph{normal form} of $g$ is a $S$-reduced sequence  $( s_{1}, \cdots, s_{n})$
such that  $g=s_1 \cdots s_n$ (If $g \in B$, a normal form of
  $g$ is $g$).  If $x= s_1\cdots s_n$, where $(s_1, \cdots, s_n)$ is  $S$-reduced, we say that $x$ is written in normal form.
\end{defi}

\begin{prop} \label{prop1}Let $(G, \ell )$ be a
valuated group with normal forms. Let ${x,y \in G}$, ${x=s_1 \cdots s_n}$ and $y=t_1 \cdots t_m$ in normal forms. If $\ell(s_nt_1) =2$, then the sequence $(s_1, \cdots, s_n, t_1, \cdots, t_m)$ is a normal form of $xy$. If ${\ell(s_nt_1)=1}$, then $(s_1, \cdots, (s_n t_1), \cdots, t_m)$ is a normal form of $xy$.
\end{prop}

\proof

Clearly if  $\ell(s_nt_1) =2$, then the sequence $(s_1, \cdots, s_n, t_1, \cdots, t_m)$ is $S$-reduced,  and therefore it is a normal form of $xy$. 

Suppose that  ${\ell(s_nt_1)=1}$. Then it is sufficient to show that $\ell(s_{n-1}s_nt_1)=\ell (s_nt_1t_2)=2$.

Suppose towards a contradiction that  $\ell(s_{n-1}s_nt_1)<2$. Then
$$
c(s_{n-1}, (s_nt_1)^{-1}) \geq \frac{1}{2},
$$
and since $c(s_{n-1}, s_n^{-1})=0$, we find by axiom $\mathcal{ A}_3$ that $$c(s_{n-1}, s_nt)=c(s_n, t_1^{-1})=0, $$ and thus $\ell(s_n)+\ell(s_nt_1)-\ell(t_1)=0$, which is clearly a contradiction.

Suppose now towards a contradiction that  $\ell(s_nt_1t_2)<2$. Then
$$
c(s_{n-1}s_n, t_1^{-1}) \geq \frac{1}{2},
$$
and since $c(t_1, t_2^{-1})=0$, we find by axiom $\mathcal{ A}_3$ that $$c(t_1, s_nt_1)=c(t_1, t_2^{-1})=0,$$ and thus $\ell(t_1)+\ell(s_nt_1)-\ell(s_n)=0$, which is clearly a contradiction.\qed

\begin{defi} Let  $(G, \ell)$ be a valuated group with normal forms.

$(1)$  An element $g \in G$ is said \emph{cyclically reduced}, abbreviated c.r., if $g \in S$ or $g=s_1 \cdots s_n$ in normal form ($n \geq 2$) and $\ell(s_ns_1)=2$. This definition does not depend on the particular choice of the normal form of $g$, as we see that if $\ell(g) \geq 2$ then $g$ is c.r. if and only if $\ell(g^2)=2 \ell(g)$.

$(2)$ An element $g \in G$ is said \emph{weakly cyclically reduced}, abbreviated w.c.r., if $g \in S$ or $g=s_1 \cdots s_n$ in normal form ($n \geq 2$) and $\ell(s_ns_1)\neq 0$. As before, this definition does not depend on the particular choice of the normal form of $g$, as we see that if $\ell(g) \geq 2$ then $g$ is w.c.r. if and only if $\ell(g^2) \geq 2 \ell(g)-1$.
\end{defi}

\begin{lem}\label{lem-decomp}  Let  $(G, \ell)$ be a valuated group with normal forms. Then every element $g$ of $G$ is conjugate to a c.r. element. Furthermore if $g \in N$ then,  there exist $x,y \in G$ such that $x$ is c.r., $g=y^{-1}xy$,  $\ell(g)=2 \ell(y)+\ell(x)$ and  $x \in N \cap S$.
\end{lem}

\proof

The proof is by induction on $\ell(g)$. The result is clear when $\ell(g) \leq 1$.

 Let $g=s_1 \cdots s_n$ in normal form, with $n \geq 2$. If $\ell(s_1s_2)=2$, then $g$ is c.r. and there is no thing to prove. If $\ell(s_ns_1)=1$, then  $g=s_n^{-1}(s_ns_1s_2 \cdots s_{n-1})s_n$. Now the sequence $(s_{n-1}, s_ns_1)$ is $S$-reduced and therefore $s_ns_1s_2 \cdots s_{n-1}$ is a c.r. element. If $\ell(s_ns_1)=0$, then $s_1=s_n^{-1}h$ for some $h \in B$ and thus $g=s_n^{-1}(hs_2 \cdots s_{n-1})s_n$. But $\ell(hs_2 \cdots s_{n-1})=n-1$, and by induction it is conjugate to a c.r. element, and the same thing holds also for $g$.

Suppose now that $g \in N$. The proof is by induction on $\ell(g)$. If $\ell(g) \leq 1$ the result is clear. Let $g=s_1 \cdots s_n$ in normal form with $n \geq 2$. Since $g \in N$ we have $\ell(g^2) \leq \ell(g)$ and therefore we get $\ell(s_ns_1)=0$; as otherwise $\ell(g^2)=2 \ell(g)-1>\ell(g)$. Then $s_n=h s_1^{-1}$ for some $h \in B$.

If $n=2$, then $g=s_1s_2=s_1 h s_1^{-1}$, and thus we have the desired conclusion.

For $n \geq 3$, we have  $g=s_1 (s_2 \cdots s_{n-1}h) s_{1}^{-1}$ for some $h \in B$. Let $g'=s_2 \cdots s_{n-1}h$. Since $\ell(g)=2+\ell(g')$, the conclusion follows by induction if we  show that  $g' \in N$.  But a simple count with normal forms shows that $\ell(g^2)=\ell(g'^2)+2 \leq \ell(g)=2+\ell(g')$, and thus $g' \in N$ as desired. \qed

\section{Subgroup Theorem for valuated groups}

The subject of this section is to prove the following theorem.
\begin{theo} \label{princip-theo}\emph{(Grushko-Neumann version  for valuated
groups)} Let $(G,\ell)$ be a valuated group with normal forms. Let $K$ be a
subgroup of $G$ such that all conjugates of  $K$ intersect
$B$ trivially. Then $K=F*_{i}G_{i}$, where  $F$ is a free group and  ${G_{i}=K \cap N(x_i)^{a_i}}$
for some $a_i\in G$ and $x_i\in S$ such that $\ell(x_i^2) \leq 1$.
\end{theo}

We introduce the following axiom:
$$
\mathcal{A}_{5}^{*}.\: \hbox{ If } \: c(x,y)+c(x^{-1},y^{-1})>\ell( x)
=\ell( y) ,\: \hbox{ then }\: \: xy^{-1}\in B^G,
$$
where  $B^G$ is the set of all conjugates of $B$.

We show first the following.
\begin{prop}
A valuated group with normal forms satisfies the axiom $\mathcal{A}_{5}^{*}$.
\end{prop}

It follows in particular, that a free product with amalgamation or an HNN-extension satisfies $\mathcal{A}_{5}^{*}$.

\proof

Let $x,y \in G$ such that
$$
c(x,y)+c(x^{-1},y^{-1})>\ell( x)
=\ell( y).
$$
Then, after simplifications, we find
$$
\ell( xy^{-1}) +\ell(
x^{-1}y) <2\ell(x).   \leqno (1)
$$

Clearly the case  $\ell( x) =\ell( y) =0$ is impossible.

Let us treat the case  $\ell( x) =\ell( y)
=1$.   Then, by (1),  $\ell( xy^{-1}) +\ell( x^{-1}y) <2$. Hence $\ell(
xy^{-1}) =0$ or $\ell( x^{-1}y) =0$.

If $\ell( xy^{-1}) =0$, then  $xy^{-1} \in B$ and we are done. If $\ell( x^{-1}y) =0$, then  $x^{-1}y=b$ for some  $b\in B$, and thus
$xy^{-1}=xb^{-1}x^{-1} \in B^G$. This ends the proof in the case   $\ell( x) =\ell( y)
=1$.

Suppose now that  $\ell( x) =\ell( y) =n\geq 2$. By (1) we have
$$
\ell( xy^{-1}) +\ell(
x^{-1}y) <2n,  \leqno (2)
$$
and in particular  $\ell( xy^{-1}) <2n$.

Let $x=s_{1}\cdots s_{n}$ and $y=t_{1}\cdots t_{n}$  in normal forms.  Since  $$\ell( xy^{-1}) <2n=\ell(x)+\ell(y),$$ we get   $\ell( s_{n}t_{n}^{-1}) \leq 1$.

If $\ell( s_{n}t_{n}^{-1}) = 1$, then $\ell(xy^{-1})=2n-1$, and by (2) we get
 $\ell( x^{-1}y)=0$.  Then  $x^{-1}y=b$ for some  $b\in B$, and thus
$xy^{-1}=xb^{-1}x^{-1} \in B^G$.

So we suppose that $\ell( s_{n}t_{n}^{-1}) =0$.

Let $i$ be the smallest  natural number such that $\ell( s_{i}\cdots s_{n}t_{n}^{-1}\cdots t_{i}^{-1}) =0$. Then we can write:
$$
x=ab,~~\ell( x) =\ell( a )+\ell( b), \hbox{ where }a=s_1 \cdots s_{i-1}, \; b=s_{i} \cdots s_{n},
$$
$$
y=\alpha \beta,~~\ell( y) =\ell( \alpha ) +\ell(\beta), \hbox{ where } \alpha=t_1 \cdots t_{i-1}, \; \beta=t_{i} \cdots t_{n}.
$$

We notice that  $\ell( a) =\ell( \alpha )$ and   $\ell( b) =\ell(
\beta )$.

Let $w=a^{-1}\alpha \beta b^{-1}$.

\bigskip
\noindent \emph{Claim 1. We have $\ell( b^{-1}w) =\ell( b^{-1})
+\ell( w)$.}

\smallskip
\noindent
\emph{Proof.}
Since  $c(x,y)+c(x^{-1},y^{-1})>\ell( x) =\ell( y) $ and  $\ell( y)
-c(x,y)=c(y^{-1},xy^{-1})$  we find
$$
c(x^{-1},y^{-1})>c(y^{-1},xy^{-1}),
$$
and by  using axiom  $\mathcal{A}_{3}$,  we get
$$c(y^{-1},xy^{-1})=c(x^{-1},xy^{-1}). \leqno (3)$$

A simplification of the expression appearing in (3) gives:
$$\ell( x^{-1}yx^{-1}) =\ell( x), \hbox{ and thus }$$$$\ell(b^{-1}wa^{-1})=\ell( b^{-1}a^{-1}\alpha \beta b^{-1}a^{-1})
=\ell( ab).\leqno (4)$$

From  $\ell( \beta b^{-1}) =0$, we get
$$\ell( b^{-1}w) =\ell(
b^{-1}a^{-1}\alpha (\beta b^{-1}))=\ell(
b^{-1}a^{-1}\alpha ),\leqno (5)
$$
$$
\ell(w)=\ell(a^{-1}\alpha (\beta b^{-1}))=\ell(a^{-1} \alpha). \leqno (6)
$$

Since  $\ell( ab) =\ell( a) +\ell( b)$ we have
$$
\ell( b^{-1}a^{-1}) =\ell( b^{-1}) +\ell( a^{-1}),
$$ and since
$\ell( a) =\ell( \alpha )$, we get
$$
\ell( b^{-1}a^{-1}\alpha
)=\ell( b^{-1}) +\ell( a^{-1}\alpha).$$

Therefore, using (5) and (6), we get
$$
\ell( b^{-1}w) =\ell( b^{-1})
+\ell( w),$$
as claimed. \qed

\bigskip
\noindent \emph{Claim 2.  We have $\ell( a^{-1}\alpha ) = 0$. }

\smallskip
\noindent
\emph{Proof.}
We treat the following two cases.

\bigskip
\noindent \emph{Case 1}. $\ell( b^{-1}wa^{-1}) =\ell( b^{-1}w) +\ell( a^{-1})$.

\bigskip
Using (4)   and Claim 1 we get
$$
\ell(a)+\ell(b)=\ell( b^{-1}wa^{-1})= \ell( b^{-1}w) +\ell(
a^{-1}) =\ell( b^{-1}) +\ell( w) +\ell( a^{-1}). $$

It follows that $\ell( w)=0$ and therefore, by (6),   $\ell( a^{-1}\alpha ) =0$ as desired.

\bigskip
\noindent \emph{Case 2}.   $\ell( b^{-1}wa^{-1}) <\ell( b^{-1}w) +\ell( a^{-1})$.

\bigskip
Then  $c(b^{-1}w,a)>0$ and since  $c(b^{-1},a)=0$, by using
$A_{2}$  we find
$$
c(b^{-1},b^{-1}w)=0.
$$

A simplification of the above expression gives
$$\ell( \beta ^{-1}\alpha ^{-1}ab) =\ell(
b)+\ell(b^{-1}w).$$

Therefore, using  (6) and Claim 1 we get
$$
\ell(x^{-1}y)=\ell(y^{-1}x)=\ell( \beta ^{-1}\alpha ^{-1}a.b)=2
\ell(b)+\ell(a^{-1}\alpha).
$$
But, counting with  normal forms we get
$$
\ell(xy^{-1})=\ell(a(b \beta^{-1})
\alpha)=\ell(a)+\ell(\alpha)\leqno (8)
$$
$$
 \hbox { or } \ell(xy^{-1})=\ell(a(b \beta^{-1}) \alpha)=\ell(a)+\ell(\alpha)-1.
$$
and thus
$$
\ell(xy^{-1})+\ell(x^{-1}y)=2
\ell(b)+\ell(a^{-1}\alpha)+\ell(a)+\ell(\alpha),
$$
$$
\hbox{ or }\ell(xy^{-1})+\ell(x^{-1}y)=2
\ell(b)+\ell(a^{-1}\alpha)+\ell(a)+\ell(\alpha)-1.
$$

Since $\ell(b)=\ell(\beta)$, we get
$$
\ell(xy^{-1})+\ell(x^{-1}y)=2n+
\ell(a^{-1}\alpha),
$$
$$
\hbox{ or } \ell(xy^{-1})+\ell(x^{-1}y)=2n+
\ell(a^{-1}\alpha)-1.
$$
By (2), the first case is impossible. Again by (2), we find  $\ell(a^{-1}\alpha)-1<0$ and
finally $\ell(a^{-1} \alpha)=0$ as desired. This ends the proof of the claim. \qed

\smallskip
Therefore $\ell( a^{-1}yx^{-1}a) =\ell( a^{-1}\alpha
\beta b^{-1}) =0$, and thus $a^{-1}yx^{-1}a=b$ for some $b\in B$. Hence
$yx^{-1}\in B^G$ as desired.  \qed

\bigskip
\noindent \textbf{Proof of Theorem \ref{princip-theo}}.

We let $K$ equipped with the induced length function. Since  $G$ satisfies  $\mathcal {A}_{5}^{*}$ and  $K$ satisfies $gKg^{-1}\cap B=\{
1\}$ for any $g \in G$, we find that  $K$ satisfies  $\mathcal{A}_{5}$ and $\mathcal{A}_{1}^*$. We notice that for any $x \in K$, $N_K(x)=K \cap N_G(x)$. Therefore by Theorem
\ref{theo-decomp-lyondon} and since $\equiv$ is an equivalence relation, we have
$$K=F*_{i}G_{i}, \hbox { where }  F \hbox{ is free and } G_{i}=K\cap N_G(y_{i}),
\hbox{ with } y_{i}\in N_G\cap K.
$$

In what follows we denote by $N(x)$ (resp. $N$) the set $N_G(x)$ (resp. $N_G$).

We are going to show that $G_{i}=K\cap N(x_i)^{a_i}$ for some
 $a_i\in K$ and  $x_i\in S\cap N$.

If  $y_{i}\in S\cap N$, then there is no thing to prove. So we suppose that $y_{i}\notin
S\cap N$.  By Lemma \ref{lem-decomp} we have
$$y_{i}=s_{1}\cdots s_{n}hs_{n}^{-1}\cdots s_{1}^{-1},$$
$$
\hbox{ with }\ell( s_{j}) =1, \;\ell( y_{i}) =2n+\ell( h), \; \ell( h) \leq 1,\; h\in N.
$$

Since $(s_{1}\cdots s_{n})^{-1}G_{i}(s_{1}\cdots
s_{n})\cap B=\{ 1\} $, we have  $\ell( h) =1$.

Let us show that  $K\cap s_{1}\cdots s_{n}N(h)s_{n}^{-1}\cdots
s_{1}^{-1}=K\cap N(y_{i})$.

We first show that $K\cap s_{1}\cdots s_{n}N(h)s_{n}^{-1}\cdots
s_{1}^{-1}\subseteq K\cap N(y_{i})$.

Let  $\beta \in N(h)$ and  $y=s_{1}\cdots s_{n}\beta
s_{n}^{-1}\cdots s_{1}^{-1}\in K$. Since $\beta \equiv h$ and $\ell(h)=1$ we have $\ell(\beta)=1$. Thus $\ell(y) \leq 2n+1$.  Let us show that $\ell(y)=2n+1$.
Suppose towards a contradiction that $\ell(y) < 2n+1$. Then, using properties of normal forms, we have
$$ \ell( s_{n}\beta ) \leq 1\quad \hbox{ or }\quad \ell(
\beta s_{n}^{-1}) \leq 1
$$
and since  $\beta $$\equiv h$, we must have
$$
\ell( s_{n}h) \leq 1\quad \hbox{or}\quad \ell( hs_{n}^{-1}) \leq 1.
$$

Therefore  $\ell( y_i) <2n+\ell( h) $, contradiction with $\ell(
y_i) =2n+\ell( h)$.

Thus  $\ell(y) =2n+1=\ell( y_{i})$ as desired. But we also have
$$ \ell(
y_{i}y) =\ell( s_{1}\cdots s_{n}h\beta s_{n}^{-1}\cdots
s_{1}^{-1}) \leq \ell( y).
$$

Therefore $y_{i}\equiv y$, hence  $y\in K\cap N(x_{i})$ as desired.

We show now that $K\cap s_{1}\cdots s_{n}N(h)s_{n}^{-1} \cdots
s_{1}^{-1}\supseteq K\cap N(y_{i})$.

Let  $y\in K\cap N(y_{i})$. Then by Lemma \ref{lem-decomp}, there exists  $\beta $ such that  $$y=t_{1}\cdots t_{n}\beta
t_{n}^{-1}\cdots t_{1}^{-1}. $$

Since  $y_{i}\equiv y$, we have
$\ell( \beta ) =1$ and  $\ell( y_{i}y) \leq \ell( y_{i})$.

Therefore
$$ \ell(
hs_{n}^{-1}\cdots s_{1}^{-1}t_{1}\cdots t_{n}\beta ) \leq 1
$$
and thus  $\ell( s_{n}^{-1}\cdots s_{1}^{-1}t_{1}
\cdots t_{n}) =0$.  Therefore  $t_{1}\cdots t_{n}=s_{1}
\cdots s_{n}\gamma$ where  $\gamma \in B$. Thus
$$
y=s_{1}\cdots s_{n}\gamma \beta \gamma ^{-1}s_{n}^{-1}\cdots
s_{1}^{-1}.
$$

Since  $\gamma \beta \gamma ^{-1}\in N$ and  $\ell(
hs_{n}^{-1}\cdots s_{1}^{-1}t_{1}\cdots t_{n}\beta )=\ell(
h\gamma \beta ) =\ell( h\gamma \beta \gamma ^{-1}) \leq 1$,  we find $h\equiv \gamma \beta \gamma ^{-1}$.

Therefore  $y\in K\cap s_{1}\cdots s_{n}N(h)s_{n}^{-1}\cdots
s_{1}^{-1}$ and this ends the proof of the theorem.  \qed

\section{A special splitting  theorem}

If $B$ is a group and $f$ is an automorphism of $B$, we denote by $B(t,f)$ the HNN-extension $\<B,t | f(b)=b^t\>$. The subject of this section is to prove the next theorem.  This theorem is needed  in the proof of Theorem \ref{thm-centralizer}. 


\begin{theo} \label{decompo-theo}
Let $(G,\ell)$ be a valuated group satisfying the following axiom: 
$$
\mathcal{A}_0^*. \hbox{ If } \ell(x) \neq 0,  \hbox{ then } \ell(x^2)>\ell(x). 
$$
Then  either $G=B$ or  there exists a sequence of automorphisms  $(f_i | i  \in \lambda)$ of  $B$, and a sequence of elements $(t_i | i \in \lambda)$ of $G$ 
such that ${G=*_B B(t_i, f_i)}$. 
\end{theo}

We will use the following lemma of Hoare.

\begin{lem} \emph{\cite{Hoare}} \label{lem-Hoare} Let $(G,\ell)$ be a valuated group. Then 

\emph{(1)} If $c(x,y)+c(x^{-1}, y^{-1}) \geq \ell(x)=\ell(y)$,  then $xy^{-1} \in N$.

\emph{(2)} Let $(u_0, \cdots, u_n)$ be a sequence of $G$ and let 
$$
a_i=u_0 \cdots u_{i-1},  i=0, \cdots, n, \hbox{ and } c_i=u_{i+1}\cdots  u_n, \quad i=1, \cdots, n-1.
$$
If 
$
\ell(a_{i+1}) \geq \ell(a_i), \ell(u_i)  \hbox{ and } \ell(c_{i-1} )\geq \ell(c_i) \hbox{ for every } i=1, \cdots, n-1,
$

$
\hbox{ and  }\ell(u_0 \cdots u_n) <\ell(u_1 \cdots u_n),
$

then 
$
\ell(a_{i+1}) = \ell(a_i),  \ell(c_{i-1}) = \ell(c_i) \hbox{ for every } i=1, \cdots, n-1$. \qed
\end{lem}

\begin{defi} Let $(G,\ell)$ be a valuated group. Let $\mathcal{U} \subseteq G$. We say that $\mathcal{U}$ is \emph{weakly reduced},  if for every sequence $(u_0, \cdots, u_n)$ of $\mathcal{U} \cup \mathcal{U}^{-1}$ which satisfies $u_iu_{i+1}\neq 1$ and $\ell(u_i) \neq 0$,  then $\ell(u_0 \cdots u_n) \geq \ell(u_1 \cdots u_n)$. 
\end{defi}

The proof of the following lemma follows the general line of the proof of \cite[Theorem(page 190)]{Hoare}.

\begin{lem} Let $(G,\ell)$ be a valuated group satisfying the axiom $\mathcal {A}_0^*$.  Let $\mathcal {U}$ be a weakly reduced subset of $G$ and let $g \in G$ such that $g \not \in \mathcal{U}^{\pm 1}$. Let $\mathcal {U}_*= \mathcal {U}\cup \{g\}$. If  $\mathcal {U}_*$ is not weakly reduced and if $\ell(g) \geq \ell(u)$ for every $u \in \mathcal{U}$, then there exists a Nielsen transformation $\phi$ of $\mathcal{U}_*$ such that $\phi(u)=u$ for every $u \in \mathcal{U}$ and $\ell(\phi(g))<\ell(g)$. 
\end{lem}

\proof Let $(u_0, u_1, \cdots, u_n)$ be a sequence of  $\mathcal {U}_*^{\pm 1}$ of minimal length for which $\mathcal {U}_*$ is not weakly reduced. Let $a_i$ and $c_i$ be the sequences as defined in Lemma \ref{lem-Hoare}. Since $n$ is minimal we have $
\ell(a_{i+1}) \geq \ell(a_i), \ell(u_i)$ and $ \ell(c_{i-1} )\geq \ell(c_i)$ for every  $i=1, \cdots, n-1
$, and therefore by  Lemma \ref{lem-Hoare},  $\ell(c_{i-1}) = \ell(c_i)$  for every $ i=1, \cdots, n-1$. 

Since $\mathcal U$ is weakly reduced, there exists $i$ such that $u_i=g^{\pm 1}$. If $i$ is unique then the transformation defined by 
$$
\phi(g)=u_0 \cdots u_n, \; \phi(u)=u \hbox{ for every } u \in \mathcal {U},
$$
is a Nielsen transformation and we have 
$$
\ell(g) \geq \ell(u_n)=\ell(c_{n-1})=\ell(c_0)>\ell(u_0 \cdots u_n)=\ell(\phi(g)). 
$$

So it is sufficient to show that the set $\{i | u_i=g^{\pm 1}\}$ is reduced to a one element. Supoose towards a contradiction that there exists $0\leq i<j\leq n$ such that $u_i=u_j=g^{\pm 1}$ or $u_i=g$ and $u_j=g^{-1}$ or $u_i=g^{-1}$ and $u_j=g$.

We have 
$$
c(c_{k-1}, c_k)=\frac{1}{2}(\ell(c_{k-1})+\ell(c_k)-\ell(u_k)), \hbox{ for every }k=1, \cdots, n-1, 
$$
and since $\ell(c_k)=\ell(c_0)$ and $\ell(u_i) \geq \ell(u_k)$ we find 
$$
c(c_{k-1}, c_k) \geq (\ell(c_0)-\frac{1}{2}\ell(u_i)), \; k=0, \cdots, n-1. \leqno (1)
$$ 
 
 We also have 
 $$
 c(c_k^{-1}, u_k)=\frac{1}{2}(\ell(c_k)+\ell(u_k)-\ell(c_{k-1}))=\frac{1}{2}\ell(u_k), \; k=1, \cdots, n-1. 
 $$
 $$
 c(c_0^{-1}, u_0)=\frac{1}{2}(\ell(c_0)+\ell(u_0)-\ell(u_0c_0)), 
$$
and since $\ell(u_0c_0)<\ell(c_0)$ we find $c(c_0^{-1}, u_0) \geq \frac{1}{2}\ell(u_0)$. Therefore 
$$
 c(c_k^{-1}, u_k) \geq \frac{1}{2}\ell(u_k), \; k=0, \cdots, n-1. \leqno (2)
$$

Now we treat the following three cases. 

\smallskip
\noindent
\textbf {Case 1}.  $u_i=u_j=g^{\pm 1}$, $0 \leq i<j\leq n-1$. 

\smallskip
By applying (2) we have 
$$
c(c_i^{-1}, u_i) \geq \frac{1}{2}\ell(u_i), \hbox{ and }c(c_j^{-1}, u_j)\geq \frac{1}{2}\ell(u_j)).
$$
But since $u_i=u_j$, by axiom $\mathcal{A}_4$, we find 
$
c(c_i^{-1}, c_j^{-1}) \geq \frac{1}{2}\ell(u_i),
$
and thus using (1) we get 
$$
c(c_i, c_j)+c(c_i^{-1}, c_j^{-1}) \geq \ell(c_i)=\ell(c_j)=\ell(c_0). 
$$

Therefore by Lemma \ref {lem-Hoare}, we find $c_jc_j^{-1} \in N$ and by axiom $\mathcal{A}_0^*$ we get $\ell(c_ic_j^{-1})=0$, and thus $\ell(u_{i+1} \cdots u_j)=0$. 

Suppose that $\ell(u_{i+2} \cdots u_j)=0$. Then $\ell(u_{i+1})=\ell(u_{i+1} \cdots u_j)=0$, a contradiction with $\ell(u_{i+1})=0$. 

Therefore  $\ell(u_{i+2} \cdots u_j)\neq 0$, and thus $\ell(u_{i+1} \cdots u_j))=0<\ell(u_{i+1} \cdots u_j))\neq 0$. Hence the sequence $(u_{i+1}, \cdots, u_j)$ contradicts the minimality of the seuqnce  $(u_{0}, \cdots, u_n)$.

\smallskip
\noindent
\textbf {Case 2}.  $u_i=u_j=g^{\pm 1}$, $0 \leq i<j=n$. 

\smallskip
Since $u_i=u_n=c_{n-1}$, then 
$$
c(c_i, c_{n-1})\geq \ell(c_0)-\frac{1}{2}\ell(u_i)=\ell(c_{n-1})-\frac{1}{2}\ell(u_n)=\frac{1}{2}\ell(c_{n-1}),
$$
and 
$$
c(c_i^{-1}, c_{n-1})=c(c_i^{-1}, u_i)=\frac{1}{2}\ell(u_i)=\frac{1}{2}\ell(c_{n-1}).
$$

Applying axiom $\mathcal{A}_4$ we find 
$$
c(c_i, c_i^{-1})=\frac{1}{2}\ell(c_{n-1})=\frac{1}{2}\ell(c_i).
$$

Therefore $\ell(c_i^2)\leq \ell(c_i)$, and thus by axiom $\mathcal{A}_0^*$, $\ell(c_i)=0$. Hence  $\ell(u_{i+1} \cdots u_{j-1})=0$. The conclusion follows as in the previous case.

\smallskip
\noindent
\textbf {Case 3}.  $u_i=u_j^{-1}$. 

\smallskip
Then $0\leq i<j-1\leq n-1$. By writting $c_n=1$ we find 
$$
c(c_{j-1}^{-1}, u_j)=\frac{1}{2}(\ell(c_{j-1})+\ell(u_{j})-\ell(c_j))\geq \frac{1}{2}\ell(u_j),
$$
and 
$$
c(c_i^{-1},u_i)\geq \frac{1}{2}\ell(u_i).
$$

By axiom  $\mathcal{A}_4$ we get  $c(c_i^{-1},c_i^{-1})\geq \frac{1}{2}\ell(u_i)$. Thus, by (1) we find 
$$
c(c_i, c_{j-1})+c(c_i^{-1}, c_{j-1}^{-1}) \geq \ell(c_i)=\ell(c_{j-1})=\ell(c_0),
$$
and thus  $\ell(u_{i+1} \cdots u_{j-1})=0$. The conclusion follows also as in the previous cases. \qed

The following lemma is a simple application of Zermilo theorem and the proof is left to the reader. 


\begin{lem} \label{lem-order}
Let $G$ be a group equipped with an integer length function $\ell : G \rightarrow \mathbb N$. Then there exists a well ordering $\prec$ of $G$ such that for every $x,y \in G$, if $\ell(x)<\ell(y)$, then $x\prec y$. \qed
\end{lem}

\noindent
\textbf{Proof of Theorem \ref{decompo-theo}}.$\;$

\smallskip
We may assume that $B \neq G$. Let $\prec$ be a well ordering of $G$ satisfying the conclusion of Lemma \ref{lem-order}. For every $g \in G$, we let $G_g$ to be the subgroup generated by the set $\{x \in G| x \prec g\}$.  We let 
$$
U=\{g \in G | g \not \in G_g\}, \quad U'=\{g \in U| \ell(g) \neq 0\}.
$$

\smallskip
\noindent
\emph {Claim 1.  $U$ generates $G$. }

\smallskip
\noindent
\emph {Proof}. 
Let $H$ be the subgroup generated by $U$ and suppose towards a contradiction that $G \neq H$. Let $a$ be the smallest element of $G$ which is not in $H$. Then every element $b \prec a$ is in $H$ and thus $a \not \in G_a$. Hence $a \in U$, a contradiction. \qed

\smallskip
\noindent
\emph {Claim 2.  $U$ is weakly reduced. }

\smallskip
\noindent
\emph {Proof}.  For every $x \in U$, we let $U_x=\{y \in U | y \prec x\}$. We show that if $U_x$ is weakly reduced then $U_x \cup \{x\}$ is weakly reduced. 

Suppose towards a contradiction that for some $x \in U$, $U_x$ is weakly reduced and $U_x \cup \{x\}$ is not weakly reduced. We see that $x \in U_x^{\pm 1}$ and for every $y \in U_x$, $\ell(x) \geq \ell(y)$. By Lemma \ref{lem-Hoare}, there exists a Nielsen tranformation $\phi$ such that $\phi(u)=u$ for any $u \in U_x$ and $\ell(\phi(x))<\ell(x)$.  Therefore $\phi(x) \in G_x$ and since $U_x \subseteq G_x$, we find $x \in G_x$; wich is a contradiction with $x \in U$. 

Thus  if $U_x$ is weakly reduced then $U_x \cup \{x\}$ is weakly reduced as required. Hence by induction on the well ordering $\prec$, $U$ is weakly reduced. \qed

\smallskip
\noindent
\emph {Claim 3. B is a normal subgroup of $G$. }

\smallskip
\noindent
\emph {Proof}.  Let $b \in B$ and $x \in G$. We have 
$$
c(bx, x^{-1}bx) =\frac{1}{2}(\ell(bx)+\ell( x^{-1}bx)-\ell(x))=\frac{1}{2}\ell( x^{-1}bx),
$$
and similarly, 
$$
c(bx, (x^{-1}bx)^{-1})=\frac{1}{2}\ell( x^{-1}bx),
$$
and by using axiom $\mathcal{A}_4$ we find 
$$
c(x^{-1}bx, (x^{-1}bx)^{-1}) \geq \frac{1}{2}\ell( x^{-1}bx). 
$$

Therefore $\ell(x^{-1}bx) \geq \ell((x^{-1}bx)^2)$ and thus by axiom $\mathcal{A}_0^*$ we get $x^{-1}bx \in B$. 

\qed

\smallskip
\noindent
\emph {Claim 4.  Let $(h_0, \cdots, h_n)$ be a sequence of $B$ and $(u_0, \cdots, u_n)$ be a sequence of $U'$. If for every $0 \leq i  \leq n-1$, $u_iu_{i+1} \neq 1$, then $\ell(h_0u_0 \cdots h_nu_n) \geq \ell(h_1u_1 \cdots h_nu_n)$.  }

\smallskip
\noindent
\emph {Proof}.
By Caim 3, $B$ is a normal subgroup of $G$ and thus $$\ell(h_0u_0\cdots h_nu_nu_n^{-1}u_{n-1}^{-1}\cdots u_0^{-1})=0, $$
and therefore, 
$$
\ell(h_0u_0 \cdots h_nu_n)=\ell(u_0\cdots u_n),
$$ 
Similarly we have  $
\ell(h_1u_1 \cdots h_nu_n)=\ell(u_1\cdots u_n).
$

Now the conclusion follows from the fact that $U$ is weakly reduced. \qed

For each $u \in U'$,  we let $F(u)=\<B,u\>$.

\smallskip
\noindent
\emph {Claim 5.  Let $f_u(x)=x^u$, restricted to $B$. Then $F(u)=B(u,f_u)$. }

\smallskip
\noindent
\emph {Proof}. It is sufficient to show that if $(h_0, \cdots, h_n)$ is a sequence of $B$ and $(u_1, \cdots, u_n)$ is a sequence of $\{u,u^{-1}\}$ such that for every  $0 \leq i  \leq n-1$, $u_iu_{i+1} \neq 1$, then $h_0u_0 \cdots h_nu_n \neq 1$. But this is a consequence of Claim 4. 

\qed

\smallskip
\noindent
\emph {Claim 6.  $G$ is the free product of $B(u,f_u)$, for $u \in U'$,  amalgamating $B$.  }

\smallskip
\noindent
\emph {Proof}. First of all, we need to show that $F(u) \cap F(v) =B$,  for $u,v \in U'$, $u \neq v$.  Let $g \in F(u) \cap F(v)$ and suppose that $\ell(g)>0$. Let $g=h_0u^{\varepsilon_1} \cdots h_{n-1}u^{\varepsilon_n}h_n$ (resp.  $g=h'_0v^{\varepsilon'_1} \cdots h'_{n-1}v^{\varepsilon'_m}h'_m$) in normal form relatively to the HNN-structure of $F(u)$ (resp. $ F(v)$).  Then
$$
h_0u^{\varepsilon_1} \cdots h_{n-1}u^{\varepsilon_n}h_n h'^{-1}_mv^{-\varepsilon'_m} \cdots v^{-\varepsilon'_1}h'^{-1}_0=1,
$$
which is a contradiction with Claim 4.  Thus  $F(u) \cap F(v) =B$ as desired.  

Now the same argument, by Calim 4, shows that $G$ is the free product  of $B(u,f_u)$, for $u \in U'$,  amalgamating $B$. \qed

\section{Conjugacy in valuated groups with normal forms}

The subject of this section is to prove the next theorem which is a generalization of \cite[Theorem 6.1]{Eric-Houcine} to valuated groups with normal forms.
\begin{theo}\label{ConjCriterion}
Let $(G,\ell)$ be a valuated group with normal forms. Let $x,y,z$
in $G$ such that $\ell(y) \geq 2$,  $y^{x}=z$. Suppose that $y$ is c.r. or w.c.r as well as $z$. Then there exist $a$, $b$ in $G$, and $n$
and $m$ in $\mathbb{Z}$, with $n \geq 1$, such that $y=(ab)^{n}$,
$z=(ba)^{n}$,  $x=a(ba)^{m}$,  and such that:

{\rm(i)} $ab$ and $ba$ are in reduced form whenever $y$ and $z$
are c.r.

{\rm(ii)} $ab$ and $ba$ are in semi-reduced form whenever $y$ and
$z$ are w.c.r.

{\rm(iii)} $ab$ is in reduced form and $ba$ is in semi-reduced
form whenever $y$ is w.c.r. and $z$ is c.r.
\end{theo}

\bigskip
The rest of this section is devoted to the proof of Theorem
\ref{ConjCriterion}, so we adopt all the assumptions and the
notation of the statement of that theorem for the rest of this
section.  The proof is in fact analogous to  the one of \cite[Theorem 6.1]{Eric-Houcine} by taking car here of the existence of elements of length 0. We shall reproduce the proof  with the necessary modifications.

We first treat the case $\ell(x)=0$.  By  putting  $b=x^{-1},a=yx$ we find $y=ab, z=ba$ and $x=b^{-1}=a(ba)^{-1}$. Since $\ell(x)=0$ we see that $ab$ and $ba$ are in reduced form.

Now, we shall treat the three cases (i)--(iii) separately. 

\medskip
\noindent \textbf{Case (i): $y$ and $z$ c.r.}

Let $x=x_1\cdots x_p$, $y=y_1\cdots y_n$, and $z=z_1\cdots z_m$ in
normal forms. We are going to prove the theorem by induction on
$p=\ell(x)$.

We first treat the case  $p=1$. Then $x^{-1}y_1\cdots y_nx=z$.
Since $y$ and $z$ are c.r. and $\ell(y) \geq 2$, we have
$\ell(x^{-1}y_1)=0$ or $\ell(y_nx)=0$. If $\ell(x^{-1}y_1)=0$, then, $x=y_1
\gamma$ for some $\gamma \in B$. Since $y$ is c.r.,
$\ell(y_ny_1)=2$, thus $\ell(y_nx)=\ell(y_ny_1\gamma)=2$ and $z=(\gamma^{-1}
y_2\cdots y_n)\cdot x$ is in reduced form. By putting
$b=\gamma^{-1} y_2\cdots y_n$ and $a=x$, we have $y=ab$, $z=ba$
and $x=a$.

Now, if $\ell(y_nx)=0$, then $x=y_n^{-1} \gamma$ for some $\gamma \in
B$.  Since $y$ is c.r., $\ell(y_ny_1)=2$, thus
$\ell(x^{-1}y_1)=\ell(\gamma^{-1}y_ny_1)=2$, then, $z=x^{-1}\cdot
(y_1\cdots y_{n-1})$ is in reduced form. By putting $a=y_1\cdots
y_{n-1}\gamma$ and $b=x^{-1}$, we have $y=ab$, $z=ba$, and
$x=b^{-1}=a(ba)^{-1}$.

We pass from $p$ to $p+1$ as follows. We have
$$x^{-1}_{p+1}\cdots  x^{-1}_1y_1\cdots y_nx_1\cdots x_{p+1}=z.$$
Since $y$ and
 $z$ are c.r. and $\ell(y) \geq 2$, we have $\ell(x_1^{-1}y_1)=0$ or $\ell(y_nx_1)=0$. We first treat the case $\ell(x_1^{-1}y_1)=0$.

\medskip
\noindent \textbf{Case (1):} $\ell(x_1^{-1}y_1)=0$.

Then $x_1^{-1}y_1=\gamma$ for some $\gamma \in B$. Then we have
$$x^{-1}_{p+1}\cdots  x^{-1}_2\gamma y_2\cdots y_ny_1\gamma^{-1}x_2\cdots x_{p+1}=z.$$
Put $x'=x_2\cdots x_{p+1}$ and $y'=\gamma y_2\cdots y_ny_1
\gamma^{-1}$. Then $y'$ is c.r. and $\ell(y') \geq 2$. By induction
there exist $a_1$, $b_1$, $\alpha \geq 1$ and $\beta$ such that
$y'=(a_1b_1)^{\alpha}$, $z=(b_1a_1)^{\alpha }$, and
$x'=a_1(b_1a_1)^{\beta}$.

\medskip
\noindent \textbf{Subcase (1-a):} $\ell(a_1)=0$ or $\ell(b_1)=0$.

Then $y'=C^{\alpha}, z=(\delta^{-1}C\delta)^{\alpha}$ and $x'=\delta(\delta^{-1}C\delta)^{s}$ for some $s\in \mathbb{Z}$,
where $C=a_1b_1$, and $\delta=a_1$ whenever $\ell(a_1)=0$ and $\delta=b_1^{-1}$ whenever $\ell(b_1)=0$.

Since $y'$ is c.r., $C$ is c.r. Thus we can write
$C=C'(y_1\gamma^{-1})$ in reduced form for some $C'$, and
$(y_1\gamma^{-1})C'$ is also in reduced form. Put
$a=y_1\gamma^{-1}\delta$ and $b=\delta^{-1}C'$. Then
$$y=y_1\gamma^{-1}y'\gamma y_1^{-1}=y_1\gamma^{-1}(C'y_1\gamma^{-1})^{\alpha}\gamma y_1^{-1}
=$$$$y_1\gamma^{-1}(C'y_1\gamma^{-1})^{\alpha-1}C'=
(y_1\gamma^{-1}C')^{\alpha}=(ab)^{\alpha},$$
$$z=(\delta^{-1}C\delta)^{\alpha}=(\delta^{-1}C'y_1\gamma^{-1}\delta)^{\alpha}=(ba)^{\alpha}$$
$$x=x_1x'=y_1\gamma^{-1}\delta(\delta^{-1}C'y_1\gamma^{-1}\delta)^s=a(ba)^s.$$

\medskip
\noindent \textbf{Subcase (1-b):} $\ell(a_1)\neq 0$ or $\ell(b_1)\neq 0$.

 Since  $\ell(b_1)\neq 0$ we can write
 $b_1=B'(y_1\gamma^{-1})$ in reduced form for some $B'$. Put $a=y_1\gamma^{-1}a_1$ and $b=B'$.
 Then $$y=y_1\gamma^{-1}y'\gamma y_1^{-1}=y_1\gamma^{-1}(a_1B'y_1\gamma^{-1})^{\alpha
 }y_1\gamma^{-1}=$$$$
 y_1\gamma^{-1}(a_1B'y_1\gamma^{-1})^{\alpha -1}a_1B'=(y_1\gamma^{-1}a_1B')^\alpha= (ab)^\alpha,$$

 $$z=(b_1b_1)^\alpha =(B'y_1\gamma^{-1}a_1)^\alpha=(ba)^\alpha, $$  where $ab$ and $ba$ are in reduced
 forms, and
 $$x=x_1x'=y_1\gamma^{-1}a_1(b_1a_1)^\beta =y_1\gamma^{-1}a_1(B'y_1\gamma^{-1}a_1)^\beta =a(ba)^\beta. $$

\noindent \textbf{Case (2):} $\ell(y_nx_1)=0$.

By taking inverses we get
$$x^{-1}_{p+1}\cdots  x^{-1}_1y_n^{-1}\cdots y_1^{-1}x_1x_2\cdots x_{p+1}=z^{-1}.$$
Therefore, by case (1), there exist $a_1$,$b_1$, $\alpha$ and
$\beta$ such that $y^{-1}=(a_1b_1)^{\alpha}$,
$z^{-1}=(b_1a_1)^{\alpha}$ and $x =a_1(b_1a_1)^{\beta}$.  Now, by
taking $a=b_1^{-1}$ and $b=a_1^{-1}$, we have $y=(ab)^\alpha$,
$z=(ba)^\alpha$ and  $$x=a_1(b_1a_1)^{\beta
}=b^{-1}(a^{-1}b^{-1})^{\beta}=aa^{-1}b^{-1}(a^{-1}b^{-1})^{\beta}=a(ba)^{-\beta-1}.$$

\medskip
\noindent \textbf{Case (ii): $y$ and $z$ w.c.r.}

Since $y$ and $z$ are w.c.r., we have $\ell(y),\ell(z) \geq 3$.  Let
$x=x_1\cdots x_p$, $y=y_1\cdots y_n$, and $z=z_1\cdots z_m$ in
normal forms. Let $y'=y_1^{-1}yy_1=y_2\cdots (y_ny_1)$ and
$z'={z_1}^{-1}zz_1=z_2\cdots (z_m z_1)$. Then $y'$ and $z'$ are
c.r. and ${z_1}^{-1}x^{-1}y_1(y')y_1^{-1}xz_1=z'$. Put
$x'=y_1^{-1}xz_1$. Then ${x'}^{-1}y'x'=z'$, $\ell(y') \geq 2$, and by
the previous case there exist $a_1$, $b_1$, $\alpha$, and $\beta$
such that $y'=(a_1b_1)^{\alpha}$, $z'=(b_1a_1)^{\alpha}$ and
$x'=a_1(b_1a_1)^{\beta}$.
\bigskip

\noindent \textbf{Case (1):} $\ell(a_1)=0$ or $\ell(b_1)=0$.

Then $y'=C^{\alpha}, z'=(\delta^{-1}C\delta)^{\alpha}$ and $x'=C^{s}\delta$ for some $s\in \mathbb{Z}$,
where $C=a_1b_1$, and $\delta=a_1$ whenever $\ell(a_1)=0$ and $\delta=b_1^{-1}$ whenever $\ell(b_1)=0$.

Since $y'$ is c.r., $C$ is c.r. Thus we can write $C=C'(y_ny_1)$
in reduced form, for some $C'$. Now since  $y_2\cdots
y_{n-1}(y_ny_1)=\delta z_2\cdots z_{m-1}(z_mz_1)\delta^{-1}$, we have $y_ny_1=\gamma
z_mz_1\delta^{-1}$ and that $y_2\cdots y_{n-1}=\delta z_2\cdots z_{m-1}\gamma^{-1}$
for some $\gamma \in B$.

Put $a=y_1C'\gamma z_m$ and $b=z_1\delta^{-1}y_1^{-1}$. Then
$$y=y_1C^{\alpha -1}C'y_n=y_1(C'y_ny_1)^{\alpha -1}C'y_n=
(y_1C'y_n)^{\alpha}=$$$$(y_1C'\gamma
z_mz_1\delta^{-1}y_1^{-1})^{\alpha}=(ab)^{\alpha}.$$

We also have
$$z=z_1(z_2\cdots z_{m-1}z_mz_1)z_1^{-1}=z_1\delta^{-1}C^{\alpha}\delta z_1^{-1}=z_1(\delta^{-1}C'y_ny_1\delta)^{\alpha}z_1^{-1}=$$
$$z_1(\delta^{-1} C'\gamma z_mz_1)^{\alpha}z_1^{-1}
=(z_1\delta^{-1} C'\gamma z_m)^{\alpha}=(z_1\delta^{-1}y_1^{-1}y_1C'\gamma
z_m)^{\alpha}=(ba)^{\alpha}.$$

We see that $y=ab$ and $z=ba$ are in semi-reduced forms. We have
$$x=y_1x'z_1^{-1}=y_1C^{s}\delta z_1^{-1}.$$
\indent If $s\geq 0$, then
$$x=y_1x'z_1^{-1}=y_1(C'y_ny_1)^{s}\delta z_1^{-1}=(y_1C'y_n)^{s}y_1\delta z_1^{-1}=$$$$(y_1C'\gamma z_mz_1\delta^{-1}y_1^{-1})^{s}y_1\delta z_1^{-1}=
(ab)^{s}b^{-1}=a(ba)^{s-1}$$

If $s<0$, then
$$x=y_1x'z_1^{-1}=y_1(y_1^{-1}y_n^{-1}C'^{-1})^{-s}\delta z_1^{-1}=(y_n^{-1}C'^{-1}y_1^{-1})^{-s}y_1\delta z_1^{-1}=$$
$$(y_1C'\gamma z_mz_1\delta^{-1}y_1^{-1})^{s}y_1\delta z_1^{-1}=
(ab)^{s}b^{-1}=a(ba)^{s-1}$$

\noindent \textbf{Case (2):} $\ell(a_1)\neq 0$ and $\ell(b_1)\neq 0$.

Since  $\ell(a_1)\neq 0$ and $\ell(b_1)\neq 0$, we can write
 $b_1=B'(y_ny_1)$  and $a_1=A'(z_mz_1)$ in reduced forms for some $B'$ and $A'$. Put
 $a=y_1A'z_m$ and $b=z_1B'y_n$.
 Then
$$y=y_1y'y_1^{-1}=y_1(a_1b_1)^{\alpha}y_1^{-1}=y_1(a_1B'y_ny_1)^{\alpha}y_1^{-1}=$$
$$y_1(a_1B'y_ny_1)^{\alpha -1}a_1B'y_n=(y_1a_1B'y_n)^{\alpha}=$$
$$(y_1A'z_mz_1B'y_n)^{\alpha }=(ab)^\alpha$$
and
$$z=z_1z'z_1^{-1}=z_1(b_1a_1)^{\alpha}z_1^{-1}=z_1(b_1A'z_mz_1)^{\alpha}{z_1}^{-1}=$$
$$(z_1b_1A'z_m)^{\alpha}=(z_1B'y_ny_1A'z_m)^{\alpha}=(ba)^\alpha,$$
and we see that $ab$ and $ba$ are in semi-reduced forms.

 If $x'=a_1(b_1a_1)^\beta $  and $\beta \geq 0$, then
 $$x=y_1x'z_1^{-1}=y_1a_1(b_1a_1)^\beta z_1^{-1}=y_1A'z_mz_1(B'y_ny_1A'z_mz_1)^\beta z_1^{-1}=$$$$
 y_1A'z_m(z_1B'y_ny_1A'z_m)^\beta =a(ba)^\beta .$$

 The case $x'=a_1(b_1a_1)^\beta $  and $\beta <0$ can be treated similarly.

\medskip
\noindent \textbf{Case (iii): $y$ w.c.r. and $z$ c.r.}

Let $x=x_1\cdots x_p$, $y=y_1\cdots y_n$, and $z=z_1\cdots z_m$ in
normal forms. Let $y'=y_1^{-1}yy_1=y_2\cdots (y_ny_1)$. Then $y'$
is c.r. and $x^{-1}y_1(y')y_1^{-1}x=z$. Put $x'=y_1^{-1}x$. Then
${x'}^{-1}y'x'=z$ and by case (i) there exist $a_1$, $b_1$,
$\alpha$, and $\beta$ such that $y'=(a_1b_1)^{\alpha}$,
$z'=(b_1a_1)^{\alpha}$ and $x'=a_1(b_1a_1)^{\beta}$. Then we
consider the case  $\ell(a_1)=0$ or $\ell(b_1)=0$, and the case $\ell(a_1)\neq 0$ and $\ell(b_1)\neq 0$. These two cases can be treated as the corresponding
subcases (1-a) and (1-b) of case (i), taking care here of the fact
that the corresponding elements $a$ and $b$ satisfy the following
condition: $ab$ is in reduced form and $ba$ is in semi-reduced
form.

This completes the proof of Theorem \ref{ConjCriterion} in all
cases. \qed

\section{Centralizer in valuated groups}

This section is devoted to study some properties of centralizers in valuated groups with normal forms.  The main subject is to show the following theorem. 


\begin{theo}\label{thm-centralizer} Let $(G,\ell)$ be a valuated group with normal forms and let $g \in G$ be  a c.r. element of length greater than 2. Then there exists a c.r. element $s$ such that $C_G(g)= \<s\> \times (B \cap C_G(g))$. 
\end{theo}

The following lemma is a detailed version of \cite[Lemma 4.9]{hur1}.  

\begin{lem} \label{lemcentrelongeure1}
Let $(G,\ell)$ be a valuated group with normal forms and let $g \in S\setminus B$. Then either
$$C_G(g)=C_G(h)^x \hbox{ for some } h, x \in G, \hbox{ such that  } h \in B \hbox{ and } x \in S\setminus B,\hbox{ or } \leqno (1)$$
$$ C_G(g)\subseteq S \cap N,\hbox{ or} \leqno (2)$$
$$C_G(g)= \<g\> \times (B \cap C_G(g)) \hbox{ and } \ell(g^2)=2. \leqno (3)$$
\end{lem}
\proof  We suppose that (1) is  not true and we show (2) or (3). We show the following claims.

\smallskip
\noindent
\emph{Claim 1. If $a, b \in S\setminus B$ and $\ell(a^{-1}ba)=1$, then $\ell(a^{-1}b)=0$ or $\ell(ba)=0$ or $\ell(a^{-1}b)=\ell(a^{-1}b)=1$}.

 \noindent
\emph{Proof.}  Suppose  that $\ell(a^{-1}b)\neq 0$ and  $\ell(ba)\neq 0$. We show  that in that case we must have $\ell(a^{-1}b)=\ell(ab)=1$.
Indeed, if $\ell(a^{-1}b)=2$ then by Proposition \ref{prop1},  the sequence $(a^{-1}, ba)$ is $S$-reduced whenever $\ell(ba)=1$,
and  the sequence $(a^{-1}, b, a))$ is $S$-reduced whenever $\ell(ba)=2$. But this contradicts $\ell(a^{-1}ba)=1$.  The situation is similar if we suppose that   $\ell(ba)=2$. \qed

\smallskip
\noindent
\emph{Claim 2.  We have $C_G(g)=\<S \cap C_G(g)\>$.}

\noindent
\emph{Proof.} Let  $x \in C_G(g)$ and $x=s_1\cdots s_n$ in normal form. We prove by induction on $n$ that $x \in  \<S \cap C_G(g)\>$. The conclusion is clear for $n=1$. We have, for $n \geq 2$,
$$s_n^{-1}\cdots s_{1}^{-1}gs_1\cdots s_n=g,\leqno (1)$$
and thus  $\ell(s_1^{-1}gs_1) \leq 1$. Since $g \notin
B^x$ for every $x \in S\setminus B$, we get ${\ell(s_1^{-1}gs_1) = 1}$.

We claim that $\ell(s_1^{-1}g)=0$ or $\ell(gs_1)=0$. If it is not the case then, by Claim 1, $\ell(s_1^{-1}g)=\ell(gs_1)=1$. But in that case, by the above proposition, since the sequence $(s_1, s_2)$ is $S$-reduced and $\ell(s_1^{-1}gs_1)=1$,  the sequence $(s_1^{-1}gs_1, s_2)$ is $S$-reduced; similarly the sequence   $(s_2^{-1}, s_1^{-1}gs_1)$ is $S$-reduced. Hence the sequence $(s_2^{-1}, s_1^{-1}gs_1, s_2)$ is $S$-reduced, and thus the sequence $$(s_n^{-1}, \cdots, s_2^{-1}, s_1^{-1}gs_1, s_2, \cdots,s_n)$$ is $S$-reduced; a contradiction with $\ell(g)=1$.

Thus $\ell(s_1^{-1}g)=0$ or $\ell(gs_1)=0$ as claimed. Hence $g=s_1h$ or $g=hs_1^{-1}$ for some $h \in B$. Thus, replacing in (1), we have
$$
s_n^{-1}\cdots s_{2}^{-1}hs_1\cdots s_n=g, \hbox{ when } g=s_1h,
$$
$$
s_n^{-1}\cdots s_{2}^{-1}s_1^{-1}hs_2\cdots s_n=g, \hbox{ when } g=hs_1^{-1},
$$
which can be rewritten as
$$
s_n^{-1}\cdots s_{2}^{-1}hgh^{-1}s_2\cdots s_n=g, \hbox{ when } g=s_1h,
$$
$$
s_n^{-1}\cdots s_{2}^{-1}h^{-1}ghs_2\cdots s_n=g, \hbox{ when } g=hs_1^{-1}.
$$

By induction,
$$ h^{-1}s_2\cdots s_n\in  \<S \cap C_G(g)\>, \hbox{ or } hs_2\cdots s_n \in  \<S \cap C_G(g)\>,$$ depending on the case $g=s_1h$ or $g=hs_1^{-1}$. Therefore
$$
s_1 \cdots s_n=(s_1h)  h^{-1}s_2\cdots s_n\in  \<S \cap C_G(g)\>, \hbox{ when } g=s_1h,
$$

$$
s_1 \cdots s_n=(s_1h^{-1})  hs_2\cdots s_n\in \<S \cap C_G(g)\>,\hbox{ when } g=hs_1^{-1},
$$
and this completes the proof of the claim. \qed

\smallskip
\noindent \emph{Claim 3.}

(i) \emph{If $\ell(g^2) \leq 1$ then $C_G(g) \subseteq S \cap N$.}

(ii) \emph{If $\ell(g^2)=2$ then $C_G(g)= \<g\> \times (B \cap C_G(g))$.}

\noindent
\emph{Proof.} $\;$

(i) By Claim 2, it is sufficient to show that if $s_1, s_2 \in S \cap C_G(g)$ then ${s_1s_2 \in S \cap N}$.

Let us show first  that if $s \in S \cap C_G(g)$ then $s \in N$. Since $\ell(s^{-1}gs)=1$, by Claim 1,  we have $\ell(s^{-1}g) \leq 1$ and $\ell(gs) \leq 1$.  Therefore
$$
c(g, s^{-1}) \geq \frac{1}{2}, c(g^{-1}, s) \geq \frac{1}{2}.
$$

Since $\ell(g^2) \leq 1$ we get $c(g,g^{-1}) \geq \frac{1}{2}$. Using axiom $\mathcal{A}_3$ we find that $c(s,s^{-1}) \geq \frac{1}{2}$ and thus $\ell(s^2) \leq 1$ as desired.

By Claim 2, it is sufficient to show that if $s_1, s_2 \in S \cap C_G(g)$ then ${s_1s_2 \in S \cap N}$. Let us show now  that $s_1s_2 \in S$. Suppose that $\ell(s_1s_2)=2$; in particular we have $\ell(s_1)=\ell(s_2)=1$. Since $\ell(s_1^{-1}gs_1)=1$, by Claim 1, we get $\ell(s_1^{-1}g)=0$ or $\ell(gs_1)=0$. Therefore $gs_2= hs_1s_2$ whenever $\ell(s_1^{-1}g)=0$ and $s_2^{-1}g=hs_2^{-1}s_1^{-1}$ whenever $\ell(s_1^{-1}g)=0$, for some $h \in B$. Hence $\ell(gs_2)=2$ or $\ell(s_2^{-1}g)=2$. Since $\ell(s_2^{-1}gs_2)=1$, again by Claim 1, $\ell(s_2^{-1}g)=0$ whenever  $\ell(gs_2)=2$; and $\ell(gs_2)=0$ whenever  $\ell(s_2^{-1}g)=2$ . We conclude that
$$
c(g^{-1}, s_2) \geq \frac{1}{2} \hbox{ and }c(g, s_2^{-1})=0, \hbox{ or } \leqno (2)
$$
$$
c(g^{-1}, s_2) =0 \hbox{ and }c(g, s_2^{-1}) \geq \frac{1}{2}.  \leqno (3)
$$

Since $c(g, g^{-1}) \geq \frac{1}{2}$ and  $c(s_2, s_2^{-1}) \geq \frac{1}{2}$ we find, using axiom $\mathcal{A}_3$,
$$
c(g, s_2^{-1}) \geq \frac{1}{2} \hbox{ if } c(g^{-1}, s_2) \geq \frac{1}{2}, \hbox{ and }
$$
$$
c(g^{-1}, s_2) \geq \frac{1}{2} \hbox{ if } c(g, s_2^{-1}) \geq \frac{1}{2},
$$
a contradiction with (2) and (3). Therefore $s_1s_2 \in S$ as desired.

(ii) By Claim 2, it is sufficient to show that if $s \in S \cap C_G(g)$ then $s \in \<g, C_G(g) \cap B\>$ and $\<g\> \cap C_G(g) \cap B=1$.

Let  $s \in S \cap C_G(g)$. We claim that $\ell(s^{-1}g)=0$ or $\ell(gs)=0$. If it is not the case then, by Claim 1, $\ell(s^{-1}g)=\ell(gs)=1$. But in that case, we have
$$
c(g^{-1},s) \geq \frac{1}{2} \hbox{ and } c(g,s^{-1}) \geq \frac{1}{2},
$$
and by using axiom $\mathcal{A}_3$ we find, since $c(g,g^{-1})=0$, $c(s,s^{-1})=0$. Therefore $\ell(s^2)=2$ and thus the sequence $(s,s)$ is $S$-reduced. By the above proposition, since  $\ell(s^{-1}gs)=1$,  the sequence $(s^{-1}gs, s)$ is $S$-reduced; similarly the sequence   $(s^{-1}, s^{-1}gs)$ is $S$-reduced. Hence the sequence $(s^{-1}, s^{-1}gs, s)$ is $S$-reduced; a contradiction with $\ell(g)=1$.

Hence  $\ell(s^{-1}g)=0$ or $\ell(gs)=0$ as claimed and thus $s \in    \<g, C_G(g) \cap B\>$. Now the fact that  $\<g\> \cap C_G(g) \cap B=1$ follows from the fact that $\ell(g^p)=|p|$ for any $p \in \mathbb Z$. This ends the proof of the claim and of the lemma.
\qed

The following lemma can  be found in \cite[Lemma 4.9(ii)]{hur1}. We give a new proof of it by using Theorem \ref{ConjCriterion}. 

\begin{lem} \label{lemmcommutation}
Let $(G,\ell)$ be a valuated group with normal forms.  Let $x,y
\in G$  such that $x$ satisfies $\ell(x^2)=2\ell(x)$, and $[x,y]=1$. Then there exist $X$ in $G$, and $h_1,h_2$
 in $B$, and $m,n \in \mathbb{Z}$, such that:
 $$x=h_1X^{n}, \quad  y=h_2X^{m}, \quad  [h_1,X]=[h_2,X]=[h_1,h_2]=1,$$
 and if $\ell(x) \geq 1, \ell(y) \geq 1$ then  $\ell(y^2)=2\ell(y)$.
 \end{lem}
\proof

We prove the lemma by induction on $\ell(x)=p$.

For $p=0$.  By taking $h_1=x, h_2=1, X=y, n=0, m=1$ we find the desired conclusion. 

For $p=1$. Since $\ell(x^2)=2\ell(x)$, by Lemma \ref{lemcentrelongeure1}, we find $$C_G(x)=\<x\> \times (B \cap C_G(x))$$ and the conclusion is clear.

 We go from $p$ to $p+1$. By Theorem \ref{ConjCriterion}, there
 exist $a$, $b$ in $G$, and $n$ and $m$ in
$\mathbb{Z}$, with $n \geq 1$, such that $x=(ab)^{n}=(ba)^{n}$,
$y=a(ba)^{m}$ and $ab$ and $ba$ are in reduced forms.

We claim that $\ell(y^2)=2\ell(y)$.  Suppose first that $\ell(a)=0$ or $\ell(b)=0$.  Then, we get,  $x=C^n=\delta^{-1}C^n\delta$ and $y=C^s\delta$ for some $s\in \mathbb{Z}$,
where $C=ab$, and $\delta=a$ whenever $\ell(a)=0$ and $\delta=b^{-1}$ whenever $\ell(b)=0$. 

Since $x$ is c.r., $C$ is c.r. We have 
$$
\ell(x^2)=\ell(\delta x^2)=\ell(x\delta x)=\ell(C^n \delta C^n), 
$$
and by using Lemma, we conclude  $\ell(C\delta C)=2\ell(C)$, and thus
$$
\ell(y^2)=\ell(C^s\delta C^s)=2|s|\ell(C)=2\ell(y),
$$
as claimed. 
Now we suppose that  $\ell(a)\neq 0$ and $\ell(b)\neq 0$. Since $x$ is c.r. and $(ab)^n=(ba)^n$ , it follows that $a$ and $b$ are c.r. Therefore 
$$
\ell(y^2)=\ell(a(ba)^ma(ba)^m)=\ell(a)+m\ell(ba)+\ell(a)+m\ell(ba)=2\ell(y),
$$
whenever $m \geq 0$, and 
$$
\ell(y^2)=\ell(a(a^{-1}b^{-1})^{-m}a(a^{-1}b^{-1})^{-m})=\ell(((b^{-1}a^{-1})^{-m-1}b^{-1})^2)=2\ell(y),
$$
whenever $m \leq 0$. This ends the proof of the claim.

Let $0 \leq t \leq n-1$ in $\mathbb N$  and $k \in \mathbb Z$, such that $m=kn+t$. Let $z=a(ba)^t$. Suppose first that $t=n-1$ and $\ell(b) =0$. Then 
$$
y=a(ba)^m=a(ba)^t(ba)^{nk}=(ab)^nb^{-1}((ba)^n)^k=xb^{-1}x^k, 
$$
and since $[b,x]=1$ we find $y=x^{k+1}b^{-1}$. Therefore by taking $h_1=1, h_2=b^{-1}, X=x$ we find the desired conclusion.

Now we suppose that $0\leq t<n-1$ or $\ell(b) \neq 0$. Then, as above, $z$ is c.r. and 
$$\ell(z)=\ell(a(ba)^t) \leq \ell(a)+t\ell(ba) <
\ell(x)=n.\ell(ba).
$$ 
Since $z$ is c.r. and $[z,x]=1$, by induction,  there exist $X$ in $G$,
and $h_1,h_2$ in $B$, such that:
 $$z=a(ba)^t=h_1X^r, \hbox{ for some }r \in \mathbb{Z},$$
 $$x=(ba)^n=h_2X^s, \hbox{ for some }s \in \mathbb{Z},$$
 $$[h_1,X]=[h_2,X]=[h_1,h_2]=1.$$
 Thus we have:

 $$x=(ba)^n=h_2X^{s}$$
$$y=a(ba)^t.((ba)^n)^k=h_1X^r(h_2X^s)^k=h_1h_2^kX^{r+sk},$$
and we find the desired conclusion. This ends the proof of the lemma. 
\qed

The following is an immediate consequence of the precedent lemma. 

\begin{cor} \label{lemcentralizer}
Let $(G,\ell)$ be a valuated group with normal forms and $g \in G$
be a c.r. element such that $\ell(g) \geq 2$. Let $x \in
C_G(g)\setminus B$. Then $x$ is c.r. and either 

\emph{(1)} $x=a^{n}$, $g=a^{m}$ for some $n,m \in \mathbb{Z}^{\#}$
and for some $a \in S\setminus B$ such that $\ell(a^{2})=2$ or,

 \emph{(2)} $\ell(x) \geq 2$. \qed
\end{cor}

Now we are ready to prove Theorem \ref{thm-centralizer}. 

\smallskip
\noindent
\textbf{Proof of Theorem \ref{thm-centralizer}.} Let $C=C_G(g)$ equipped with the induced length function. Then $(C, \ell)$ is a valuated group. We claim that $(C, \ell)$ satisfies the axiom $\mathcal{A}_0^*$ of Theorem \ref{decompo-theo}. If $x \in C$, then by Corollary \ref{lemcentralizer} either $x \in B \cap C$ or $x$ is c.r. and hence $\ell(x^2) >\ell(x)$. Therefore, if $x \not \in B'=   B \cap C$, then  $\ell(x^2) >\ell(x)$; thus $(C, \ell)$ satisfies  $\mathcal{A}_0^*$ as claimed. Therefore by Theorem \ref{decompo-theo}, $C=*_{B'}B'(t_i,f_i)$, where $B'(t_i,f_i)=\<B',t_i|f(x)=x^{t_i}\>$. But since the center of $C$ contains a c.r. element, we get that $C= \<B',t_i|f(x)=x^{t_i}\>$, which can be written simply as $C=\<B',s|B'^s=B'\>$. Again, since $Z(C)$ contains a c.r. element, we find that $C=\<s\>\times B'$, as desired. \qed

We end this section with the following theorem. We give a proof of it, using only  Corollary \ref{lemcentralizer} and Theorem \ref{princip-theo}.

\begin{theo} \label{theocentralizer}
Let $(G,\ell)$ be a valuated group with normal forms. Let $g \in
G$ be a c.r. element of length greater than $2$. Then the
following properties are equivalent:

\emph{(1)} $C_G(g)\cap B=1$.

\emph{(2)} $C_G(g)$ is infinite cyclic.

\emph{(3)} $C_G(g)$ is locally cyclic.

\noindent
 In that case, if $G$ has no involutions then $C_G(g)$ is
selfnormalizing.
\end{theo}

\proof

(1)$\Rightarrow$(2) By Corollary \ref{lemcentralizer}, for every $x
\in C_G(g)$, $x \neq 1$, $x$ is c.r. and  we have:

(i) $x=a^{n}$, $g=a^{m}$ for some $n,m \in \mathbb{Z}^{\#}$
and for some $a \in S\setminus B$ such that $\ell(a^{2})=2$ or,

(ii) $\ell(x) \geq 2$.

\smallskip
 We claim that $C_G(g) \cap B^y=1$, for every $y \in G$.
 Suppose towards a contradiction that $C_G(g) \cap B^{y} \neq
 1$ for some $y \in G$ and let $z \in C_G(g) \cap B^{y}$.
 Then by the above property $z$ is c.r. and $\ell(z^{2})>\ell(z)$, a contradiction.

 Hence by Theorem \ref{princip-theo}, $C_G(g)$ is a free product. The result follows.

 (2)$\Rightarrow$(3) Obvious.

(3)$\Rightarrow$(2) We claim that $C_G(g) \cap
B^{y}=1$, for every $y \in G$.
 Suppose towards a contradiction that $C_G(g) \cap B^{y} \neq
 1$ for some $y \in G$ and let $z \in C_G(g) \cap B^{y}$.
 Then $y=b^{n}$ and $g=b^{m}$ for some $n,m \in \mathbb{Z}^{\#}$
 and for some $b$. Now since $g$ is c.r., $b$ is c.r. Thus $z$ is
 c.r., a contradiction. Hence by Theorem \ref{princip-theo}, $C_G(g)$ is a free product. The result follows.

(2) $\Rightarrow$(1) Obvious.

Let us show now that $C_G(g)$ is selfnormalizing. Let
$C_G(g)=\<s\>$, for some $s \in G$. Then we see that $s$ is c.r.
and  is not a proper power.

Let $x \in N(C_G(g))^{\#}$, then $s^{x}=s^{m}$ for some $m \in
\mathbb{Z}^{\#}$. Now since $s$ is not a proper power and c.r. we
have $m= \pm 1$. If $m=1$ we have the result. If $m=-1$ then
$x^{-2}sx^2=s$ and thus $x^{2} \in C_G(g)$. Now since $G$ has no
involution, $x^{2} \neq 1$. And since $C_G(g)\cap B=1$,
$\ell(x^{2}) \neq 0$, and by Lemma \ref{lemcentralizer},  $x^2$ is
c.r. Hence $x$ is c.r. and  $\ell(x^{2})>\ell(x)$. 

Hence, as above,  $ N(C_G(g)) \cap B^{y}=1$, for every $y
\in G$.

Thus by Theorem \ref{princip-theo}, $A=N(C_G(g))$ is a free group. But since
$N_A(C_G(g))=A$ we have $A$ is infinite cyclic. Hence $A$ is
generated by $s$ since $s$ is not a proper power.
 \qed
\section{The CSA property in valuated groups}
 If $G$ is a group, $S$ a subset of $G$ and $H$ is a subgroup, we say that $H$ is \textit{$S$-malnormal} if $H \cap H^s = 1$ for any $s \in S$, $s \neq 1$. The subject of this section is to prove the following theorem.

 \begin{theo} \label{TheoCSA} Let $(G,\ell)$ be a valuated group
 with normal forms and without involutions. Then the following properties are
 equivalent:

\noindent \emph{(1)} $G$ is a CSA$^*$-group.

 \noindent \emph{(2)} The following properties are satisfied:

 $(i)$ for every $g \in G$, $g\neq 1$, if  $~C_G(g) \subseteq S$ then $C_G(g)$ is abelian and $S$-malnormal,

$(ii)$ for every $g \in B$, $g \neq 1$, $C_G(g)$ is abelian and malnormal.
 \end{theo}

\proof

Obviously we have (1)$\Rightarrow$(2). Assume (2) and let us show  (1).
 We are going to prove that for every $g \in
G^{\#}$, $C_G(g)$ is abelian and selfnormalizing. Let us treat first the
case when $ \ell(g) \leq 1$. The case $\ell(g)=0$ follows from the
assumption $(ii)$. Suppose that  $\ell(g)=1$. By Lemma
\ref{lemcentrelongeure1}, there is three cases to consider:

\smallskip
\noindent $(a)$ $C_G(g)=C_G(h)^x \hbox{ for some } h, x \in G,
\hbox{ such that  } \ell(h)=0 \hbox{ and } \ell(x)=1, \hbox{ or}$

\noindent$ (b)$ $C_G(g) \subseteq S,\hbox{
or }$

\noindent $(c)$ $C_G(g)= \<g\> \times (B \cap C_G(g))
\hbox{ and } \ell(g^2)=2.$

\smallskip
The case $(a)$ follows form the assumption $(ii)$. Let us treat the case
$(b)$. By assumption $(i)$, $C_G(g)$ is abelian and thus we prove that
it is malnormal. We suppose also that  $g \notin B^x$ for
every $x \in G$ such that $\ell(x)=1$, for if we are in the case
$(a)$. 

Let  $x\in G$ and  $g',g'' \in C_{G}(g)^{\#}$ such that
$g'^{x}=g''$. If $x\in S$, then $x \in
C_{G}(g)$ because $C_G(g)$ is $S$-malnormal. Suppose now  that $x \notin S$,
i.e. that $x=s_1\cdots s_n$ is in normal form with $n\geq 2$. Then
$$s_n^{-1}\cdots s_1^{-1}g'
s_1\cdots s_ng''^{-1}=1.$$ This implies that
$\ell(s_1^{-1}g's_1)=0$, or $\ell(s_1^{-1}g's_1)=1$.

If $\ell(s_1^{-1}g's_1)=0$ then $g' \in C_G(\gamma)^{-s_1}$ for
some $\gamma \in B$. Now since ${g' \in
C_G(\gamma)^{-s_1}}$ and the last group is abelian and malnormal,
${g \in C_G(\gamma)^{-s_1}}$, thus ${C_G(g) \subseteq
C_G(\gamma)^{-s_1}}$. Since $C_G(\gamma)^{-s_1}$ is malnormal, $x
\in C_G(\gamma)^{-s_1}$. Finally $x \in C_G(\gamma)^{-s_1}$, as
the last group is abelian.

We claim that the case $\ell(s_1^{-1}g's_1)=1$ cannot occur. For
if $\ell(s_1^{-1}g's_1)=1$, then $g'=s_1$ and
$\ell(s_2g's_2^{-1})=1$. (If  $\ell(s_2g's_2^{-1})=0$ then we are
in the previous case). If $n=2$, then $s_2g's_2^{-1}=g''$ and $s_2
\in C_G(g)$ and $\ell(g's_1)=2$, a contradiction. If $n>2$ then
$g'=s_2$ and  $\ell(g'^2)=2$, which is also a contradiction as
$g'^2 \in C_G(g)$. This completes the proof of the case $(b)$.

We treat now the case $(c)$.

If $B \cap C_G(g) \neq 1$, then let $h \in (B \cap
C_G(g))^{\#}$. Then $C_G(g) \subseteq C_G(h)$, as $C_G(h)$ is
abelian and malnormal. We see that also in this case that $C_G(g)$
is malnormal.

If $B \cap C_G(g) = 1$,  then $C_G(g)=\<g\>$, and by Theorem
\ref{theocentralizer}, $C_G(g)$ is selfnormalizing.

Let us now treat the case when $\ell(g) \geq 2$. Then
$g=sg's^{-1}$, where $g'$ is c.r. Now if $\ell(g') \leq 1$ then,
up to conjugacy, we are in the previous case. Thus suppose that
$\ell(g') \geq 2$. Then we see as before that if $B \cap
C_G(g) \neq 1$ ,then $C_G(g)$ is abelian and selfnormalizing, and
if $B \cap C_G(g) = 1$,  then by Theorem
\ref{theocentralizer}, $C_G(g)$ is infinite cyclic and
selfnormalizing. \qed

\bibliographystyle{alpha}
\bibliography{biblio}
\end{document}